\documentclass [12pt,a4paper] {article}

\usepackage {natbib}
\usepackage {amsmath}
\usepackage {amssymb}
\usepackage {amsthm}
\usepackage {verbatim}

\DeclareMathOperator {\cl} {cl}
\DeclareMathOperator {\dcl} {dcl}
\DeclareMathOperator {\Aut} {Aut}
\DeclareMathOperator {\tp} {tp}
\DeclareMathOperator {\dom} {dom}
\newcommand {\liff} {\leftrightarrow}
\newcommand {\nforks} {\mathrel{\raise0.2ex\hbox{\ooalign{\hidewidth$\vert$\hidewidth\cr\raise-0.9ex\hbox{$\smile$}}}}}
\newcommand {\forks} {\not \nforks}
\DeclareMathOperator {\CB} {CB}
\DeclareMathOperator {\Mlt} {Mlt}

\theoremstyle {definition}
\newtheorem {definition} {Definition} [section]

\newtheorem {lemma} [definition] {Lemma}
\newtheorem {theorem} [definition] {Theorem}
\newtheorem {proposition} [definition] {Proposition}
\newtheorem {corollary} [definition] {Corollary}

\newtheorem {fact} [definition] {Fact}
\theoremstyle {remark}

\title {Constructing Quasiminimal Structures}
\author {Levon Haykazyan}

\begin {document}

\maketitle

\begin {abstract}
Quasiminimal structures play an important role in non-elementary categoricity.
In this paper we explore possibilities of constructing quasiminimal models of a
given first-order theory. We present several constructions with increasing
control of the properties of the outcome using increasingly stronger assumptions
on the theory. We also establish an upper bound on the Hanf number of the
existence of arbitrarily large quasiminimal models.
\end {abstract}

\section {Introduction}

An uncountable structure is called {\em quasiminimal} if every (first order)
definable subset is either countable or cocountable. Quasiminimal structures
carrying a homogeneous pregeometry play an important role in non-elementary
categoricity (see \cite {zilber-qme}). Several key analytic structures such as
$\mathbb C_{\exp}$ are conjectured to be quasiminimal (see \cite {zilber-pe}).

The study of quasiminimal structures from the first order perspective is
pioneered in \cite {pilovic}. They isolate a key notion of a strongly regular
type for an arbitrary theory (see Definition \ref {sreg-type}). An important
result is that assuming the generic type (i.e. the type containing formulas
defining cocountable sets) of a quasiminimal structure is definable, its unique
global heir is strongly regular. In this paper we establish the converse of
this: if there is a definable strongly regular type, then we can construct a
quasiminimal model. The definability condition holds in particular for groups.
Thus given a regular group, there is a quasiminimal group elementarily
equivalent to it. (The converse stating that the monster model of a quasiminimal
group is regular is due to \cite {pilovic}.) This reduces the existence of
non-commutative quasiminimal groups and quasiminimal fields that are not
algebraically closed to respective problems for regular groups and fields. See
the next section for a detailed discussion of this.

The techniques of \cite {pilovic} provide more. If the cardinality of a
quasiminimal structure is strictly greater than $\aleph_1$, then the global heir
of the generic type is symmetric (i.e. Morley sequence are totally
indiscernible). Thus to construct quasiminimal models of arbitrarily large
cardinalities we need additional assumptions. Here we present two constructions
of arbitrarily large quasiminimal models: one assuming the theory has definable
Skolem functions and the other assuming the theory is stable. We also use a
variant of Skolem functions to prove an upper bound on the Hanf number of having
arbitrarily large quasiminimal models.

If we strengthen the stability assumption to $\omega$-stability, we can use the
existence of prime models. We can then easily get quasiminimal models by taking
prime models over Morley sequences in the strongly regular type. Here we show
that a stronger conclusion holds: the class of all such models is quasiminimal
excellent in the sense of \cite {zilber-qme}. Since this class is clearly
uncountably categorical, its excellence was expected. However Shelah's results
on deducing excellence from categoricity require additional set theoretic
assumptions and only apply to $L_{\omega_1, \omega}$ sentences. So they can't be
used to deduce excellence in this generality.

\section {Preliminaries}

The model theoretic notation is standard. The notions left undefined are either
standard or can be found in \cite {pilovic} (or both) on which this paper is
based on. 

We work in a monster model $\mathfrak C$ of a first order theory $T$ in a
countable language. Types over $\mathfrak C$ are called global types. A global
type $\mathfrak p(x) \in S_1(\mathfrak C)$ is an ultrafilter on the boolean
algebra of definable subsets of $\mathfrak C$. If we think of sets in $\mathfrak
p$ as ``large'' and sets outside $\mathfrak p$ as ``small'', we can, in analogy
with the algebraic closure, define the closure of $A$ to be the union of all
small sets definable over $A$.  More formally define
$$\cl_{\mathfrak p}(A) = \{b \in \mathfrak C : b \not \models \mathfrak p|A\}.$$
A natural question is when is $\cl_{\mathfrak p}$ a closure operator
(monotone, idempotent, finitary operator) or better a pregeometry (closure
operator with exchange). To answer that question we need the notion of a
strongly regular type.

\begin {definition}
\label {sreg-type}
Let $A$ be a subset, $\mathfrak p(\bar x)$ be a global $A$-invariant type and
$\phi(\bar x) \in \mathfrak p$ be a formula over $A$.  The pair $(\mathfrak p,
\phi)$ is called {\em strongly regular} if for all $B \supseteq A$ and $\bar a$
satisfying $\phi$ either $\bar a \models \mathfrak p|_B$ or $\mathfrak p|_B
\vdash \mathfrak p|_{B\bar a}$. 
\end {definition}

Note that the following statements are all equivalent to $\mathfrak p|_B \vdash
\mathfrak p|_{B\bar a}$.
\begin {itemize}
\item The type $\mathfrak p|_B$ has a unique extension to a type over $B\bar a$.
\item If $\bar b \models \mathfrak p|_B$, then $\tp_{\bar x}(\bar b/B) \cup
	\tp_{\bar y}(\bar a/B)$ determines a complete type over $B$ in $\bar x
	\bar y$.
\item If $\bar b \models \mathfrak p|_B$, then $\tp(a/B) \vdash \tp(a/Bb)$.
\end {itemize}

In our case $\mathfrak p$ will be a type in one variable, $\phi(x)$ will always
be the formula $x=x$ (so we surpass $\phi$ from the notation) and we will assume
$A = \emptyset$ (this can be achieved by extending the language, if $A$ is
countable). 

\begin {fact} [\cite {pilovic}]
Let $\mathfrak p \in S_1(\mathfrak C)$ be an $\emptyset$-invariant type. Then
$\cl_{\mathfrak p}$ is a closure operator if and only if $\mathfrak p$ is
strongly regular. Moreover, $\cl_{\mathfrak p}$ is a pregeometry if in addition
all (equivalently some) Morley sequences in $\mathfrak p$ are totally
indiscernible.
\end {fact}

Here by a Morley sequence in $\mathfrak p$ we mean a sequence of $(a_i : i <
\alpha)$ satisfying $a_i \models \mathfrak p|_{\{a_j : j < i\}}$. Since
$\mathfrak p$ is $\emptyset$-invariant, a Morley sequence is always
indiscernible. Also it is easy to see that for the sequence to be totally
indiscernible it is necessary and sufficient that $a_1a_0 \equiv a_0a_1$. In
this case $\mathfrak p$ is called {\em symmetric} and it is called {\em
asymmetric} otherwise.

We can also start from the other end. Let $M$ be a structure and $\cl$ an
infinite dimensional pregeometry on it. Assume that the pregeometry is related
to the language in the following way: for every finite $B \subset M$ the set of
elements of $M$ outside $\cl(B)$ is the set of realisations of a complete type
$p_B(x)$ over $B$. Such a pregeometry is called {\em homogeneous} in \cite
{pilovic}. In the context of quasiminimal excellent classes this property is
referred as the {\em uniqueness of the generic type}. Now consider $p_{\cl}(x) =
\bigcup_{B \subset_{fin} M} p_B(x)$. It can be seen that $p_{\cl}(x) \in S_1(M)$
is a complete type over $M$. Results of \cite {pilovic} show that this type
extends to a global strongly regular type.

\begin {fact} [\cite {pilovic}]
\label {hompreg}
Let $(M, cl)$ be a homogeneous pregeometry. Then $p_{\cl}$ (defined above) is
definable and its unique global heir $\mathfrak p$ is $\emptyset$-invariant,
generically stable, symmetric strongly regular type. In particular $\cl =
\cl_{\mathfrak p}|_M$.
\end {fact}

Thus homogeneous pregeometries (such as those in quasiminimal excellent classes)
come from symmetric regular types and conversely symmetric regular types induces
homogeneous pregeometries. On the other hand we can also find strongly regular
types if we start from quasiminimal structures. 

Let $M$ be a quasiminimal structure. Then 
$$\{\phi(x) \in L(M) : \phi \text { defines a cocountable set}\}$$ 
is a complete type over $M$. We call this type the {\em generic type} of $M$.
Under suitable conditions we can also extend this type to a strongly regular
global type.

\begin {fact} [\cite {pilovic}]
\label {pilovic-quasi}
Assume that the generic type $p$ of a quasiminimal structure $M$ is
$\emptyset$-definable, then its unique global heir $\mathfrak p$ is strongly
regular. Moreover if $|M| > \aleph_1$, then the definability condition holds
(after adding parameters) and $\mathfrak p$ is symmetric and generically stable.
\end {fact}

However there are no guarantees that the strongly regular type $\mathfrak p$ is
symmetric (unless $|M| > \aleph_1$). For example the structure $(\aleph_1 \times
\mathbb Q, <)$ in the lexicographic order is quasiminimal. Indeed it has
quantifier elimination (as a dense linear order) and every formula of the form
$x < a$ defines a countable set, while every formula of the form $x > a$ defines
a cocountable set. Now the generic type is the completion of $\{x > a : a \in
\aleph_1 \times \mathbb Q\}$, so is $\emptyset$-definable. But its unique global
heir is clearly asymmetric (no indiscernible sequence is totally indiscernible
because of the order).

The definability condition holds in particular if $(M, \cdot)$ is a quasiminimal
group. The generic type is definable using the fact that a definable set $X
\subseteq M$ is cocountable if and only if $X \cdot X = M$. Indeed, if $X$ is
countable, then so is $X \cdot X$ and hence it can't be the whole of $M$.
Conversely if $X$ is cocountable, then so is $mX^{-1}$ for a given $m \in M$.
Hence $mX^{-1} \cap X$ is nonempty and therefore $m \in X \cdot X$. Thus using
the Fact \ref {pilovic-quasi}, we see that the monster model of $M$ has a
strongly regular type. Such groups are called {\em regular} in \cite {pilovic}.

One of the questions \cite {pilovic} asks is whether every regular group is
commutative. A similar question of whether every quasiminimal group is
commutative has been around for some time. By the above, given a non-commutative
quasiminimal group, its monster model is a regular non-commutative group. 
Our results imply the converse of this: if there is a regular non-commutative
group, there is a quasiminimal non-commutative group. \cite {pilovic} eludes to
a possible construction of a regular non commutative group, but to the best of
our knowledge the problem is still open.

A related question is whether every regular (or quasiminimal) field is
algebraically closed. It is known for regular fields where the strongly regular
type is symmetric and hence for quasiminimal fields of cardinality $> \aleph_1$
(see \cite {hls} and \cite {gogacz}). However the case of an asymmetric regular
field and a quasiminimal field of cardinality $\aleph_1$ remain open. Our
results together with Fact \ref {pilovic-quasi} reduce these two cases to each
other.

In the rest of the paper we are mainly concerned with constructing quasiminimal
models. From the above we know that in well behaved quasiminimal structures
(certainly for large quasiminimal structures) the generic type is strongly
regular and induces a closure operator. So a natural question is given a theory
with a strongly regular type, can we construct a quasiminimal model. Note also
that a model will be quasiminimal if and only if it satisfies the countable
closure property (i.e. closure of a countable set is countable). We present
several constructions with increasing control of properties of the outcome using
increasingly stronger assumptions on the theory.

\section {Arbitrary Theories}
\label {gensec}

In this section we describe a method for constructing a quasiminimal structures
which we then apply to construct a quasiminimal model of a theory with  a
strongly regular type. Our method is based on the original method of \cite
{morleyvaught} of constructing an $(\aleph_1, \aleph_0)$ model from a Vaughtian
pair. The method is widely known among model theorists (although not used often
outside of its original context) and expositions are available in many model
theory texts.

\begin {definition}
A {\em special pair} for a theory $T$ is a pair of models $M \prec N$ where $N$
is a proper elementary extension and there is an $\emptyset$-definable type $p
\in S_1(M)$ such that all $a \in N \setminus M$ realise $p$.
\end {definition}

Let $M \prec N$ be a special pair of models of $T$. Let $p \in S_1(M)$ be the
type of elements in $N \setminus M$ and $d_p$ be the schema defining $p$ (note
that all schemas defining $p$ are equivalent). Add a new predicate symbol $R$ to
the language and consider the theory $\hat T$ of $N$ where $R$ is interpreted as
$M$. The theory $\hat T$ in particular encodes the following
\begin {itemize}
\item The universe is a model of $T$;
\item $R$ defines a proper elementary substructure;
\item $\{\phi(x, \bar a) : \bar a \in R \land d_p x \phi(x, \bar a)\}$ is a
	complete type over $R$;
\item $\forall x \not \in R~\forall \bar y \in R (\phi(x, \bar y) \liff
	d_p x \phi(x, \bar y))$ (every element outside $R$
	realises the above type).
\end {itemize}
Thus any model of $\hat T$ provides a special pair for $T$. So if we have a
special pair for $T$, we can construct a countable one. Therefore without lose
of generality assume that $N$ is countable. Further by iteratively realising
types in $M$ and $N$ we may assume that both $M$ and $N$ are homogeneous and
realise the same types over $\emptyset$. Indeed we can construct a sequence
$$(N_0, M_0) \preccurlyeq (N_1, M_1) \preccurlyeq (N_2, M_2) \preccurlyeq ...$$
of countable models of $\hat T$ such that
\begin {itemize}
\item if $p \in S_n(T)$ is realised in $N_{3i}$, then it is realised in
	$M_{3i+1}$;
\item if $\bar a, \bar b, c \in M_{3i+1}$ and $\tp^{M_{3i+1}}(\bar a) =
	\tp^{M_{3i+1}}(\bar b)$, then there is $d \in M_{3i+2}$ such that
	$\tp^{M_{3i+2}}(\bar a, c) = \tp^{M_{3i+2}}(\bar b, d)$;
\item if $\bar a, \bar b, c \in N_{3i+2}$ and $\tp^{N_{3i+2}}(\bar a) =
	\tp^{N_{3i+2}}(\bar b)$, then there is $d \in N_{3i+3}$ such that
	$\tp^{N_{3i+3}}(\bar a, c) = \tp^{N_{3i+3}}(\bar b, d)$;
\end {itemize}
Then take $(N, M) = \bigcup_{i < \omega} (N_i, M_i)$. Now $M$ and $N$ would be
homogeneous by the second and third close and realise the same types over
$\emptyset$ by the first close. In particular $M$ and $N$ would be isomorphic.
This is a standard step in Vaught's two cardinal theorem and more details can be
found in standard textbooks such as \cite {marker}.

\begin {proposition}
If $T$ has a special pair, then it has a quasiminimal model of cardinality
$\aleph_1$.
\end {proposition}

\begin {proof}
Let $M \prec N$ be a special pair, $p$ the type of elements in $N \setminus M$ and
$d_p$ the defining scheme. By the above discussion we may assume that $M$ and
$N$ are countable, homogeneous and isomorphic. We construct an elementary
sequence $(A_\alpha : \alpha < \aleph_1)$ of structures, each isomorphic to $M$
(and hence $N$) such that for each $\beta < \alpha$ every element of $A_\alpha
\setminus A_\beta$ realises the unique heir of $p$ over $A_\beta$. The
construction is by transfinite recursion.
\begin {itemize}
\item Take $A_0 = M$.
\item If $\delta$ is a limit ordinal, then take $A_\delta = \bigcup_{\alpha <
	\delta} A_\alpha$. Note that $A_\delta$ is homogeneous, countable and
	realises the same types (over $\emptyset$) as $M$. Hence $A_\delta \cong
	M$. For $\beta < \delta$ every element in $A_\delta \setminus A_\beta$
	is in some $A_\alpha$ with $\beta < \alpha < \delta$. Therefore it
	realises the unique heir of $p$ over $A_\beta$.
\item Assume $A_\alpha$ is constructed. Let $f : M \to A_\alpha$ be an
	isomorphism. Then $f$ extends to an isomorphism $\hat f : N \to B$ for
	some structure $B$. Take $A_{\alpha+1} = B$. Then $A_\alpha \prec
	A_{\alpha+1}$ and $A_{\alpha+1} \cong N \cong M$. It remains to show
	that every $c \in A_{\alpha+1} \setminus A_\alpha$ realises the unique
	heir of $p$ over $A_\alpha$. Let $\bar a \in A_\alpha$ and $A_{\alpha+1}
	\models d_p x \phi(x, \bar a)$. Let $c' = \hat f^{-1}(c), \bar a' = \hat
	f^{-1}(\bar a)$. Note that $c' \in N \setminus M$. Since $\hat f$ is an
	isomorphism, we have $N \models d_p x \phi(x, \bar a')$. But $\hat f$
	extends $f$ and since $\bar a \in A_\alpha$, we have $\bar a' =
	f^{-1}(\bar a) \in M$. Therefor $N \models \phi(c', \bar a')$ and so
	$A_{\alpha+1} \models \phi(c, \bar a)$.
\end {itemize}
Now take $A = \bigcup_{\alpha < \aleph_1} A_\alpha$. If $\bar a \in A$, then
$\bar a \in A_\alpha$ for some countable $\alpha$. Now if $d_p x \phi(x, \bar
a)$ is satisfied, then every element in $A \setminus A_\alpha$ satisfies
$\phi(x, \bar a)$ and hence it defines a cocountable set. Otherwise $\phi(x,
\bar a)$ defines a countable set. Thus $A$ is quasiminimal. Also observe that
countability is weakly definable in $A$ in the sense of \cite {zilnotes}.
\end {proof}

Now we can apply this method to construct a quasiminimal model from a strongly
regular type.

\begin {theorem}
\label {gen-main}
Let $T$ be a theory with $\mathfrak p$ a $\emptyset$-definable strongly regular
type.  Then there is a quasiminimal model of $T$.
\end {theorem}

\begin {proof}
Let $M \prec \mathfrak C$ be a small model. Consider $N = \cl_{\mathfrak p}(M)$.
Since $\mathfrak p|_M$ must have a realisation, we have $N \subsetneq \mathfrak
C$. We claim that $N$ is an elementary substructure. We use the Tarski-Vaught
test. Assume that $\bar a \in N$ and $\phi(b, \bar a)$ holds for some $b \in
\mathfrak C$. If $b \in N$, then we are done. So assume that $b \not \in N$.
Then $b$ realises $\mathfrak p|_M$ and hence $\mathfrak p|_N$. Now since $\bar a
\in N$, we have that $\mathfrak p|_M \vdash \mathfrak p|_N$ and therefore some
$L(M)$-formula $\psi(x)$ implies $\phi(x, \bar a)$.  Since $M$ is a model,
$\psi$ must be realised in $M$ and hence in $N$. Thus $\phi$ is realised in $N$
and $N \prec \mathfrak C$. Now every element of $\mathfrak C \setminus N$
realises $\mathfrak p|_N$ and hence $N$ and $\mathfrak C$ are a special pair.
\end {proof}

We can apply this to regular groups, since by \cite {pilovic} the strongly
regular type in a group is $\emptyset$-definable.

\begin {corollary}
If $G$ is a regular group, then there is a quasiminimal group $H$ elementarily
equivalent to $G$.
\end {corollary}

Without additional assumptions we cannot construct quasiminimal models of
cardinalities larger than $\aleph_1$. Indeed by Fact \ref {pilovic-quasi} the
generic type of such a structure is symmetric. So we at least need to assume
that $\mathfrak p$ is symmetric.

\section {Theories with definable Skolem functions}

We say that a theory $T$ has {\em definable Skolem functions} if for every
$L$-formula $\phi(y, \bar x)$ there is a $\emptyset$-definable function
$f_\phi(\bar x)$ such that $T \models \forall \bar x (\exists y \phi(y, \bar x)
\to \phi(f_\phi(\bar x), \bar x))$. In this section we study symmetric strongly
regular types in theories with definable Skolem functions. So let $T$ be such
a theory and let $\mathfrak p$ be a symmetric strongly regular type. Note that
by the Fact \ref {hompreg} the type $\mathfrak p$ is $\emptyset$-definable. 

We will show that $T$ has a quasiminimal model of arbitrary uncountable
cardinality $\kappa$. Since there are definable Skolem functions, for every
subset $A$ there is the least model containing it, namely the submodel with
universe $\dcl(A)$. Thus if there is a quasiminimal model containing a Morley
sequence $A = (a_\alpha : \alpha < \kappa)$, then $\dcl(A)$ would be one.  Thus
we need to show that $\dcl(A)$ is quasiminimal. Further since $\dcl(A)$ embeds
in every model containing $A$, all we need is to construct a model $M$
containing $A$ such that for a fixed countable $A_0 \subseteq A$ we have that
$\cl_{\mathfrak p}(A_0) \cap M$ is countable. For that we will use the technique
of self-extending models originally due to \cite {vaught-self}.

The type $\mathfrak p$ is an ultrafilter on definable subsets. We can use this
ultrafilter for a definable analogue of the ultrapower construction. More
specifically given a model $M$ of $T$ we can associate a canonical extension
$M^*$ defined as follows. The elements of $M^*$ are all definable (with
parameters) functions $f : M \to M$ modulo agreeing on a large set (i.e. a
member of $\mathfrak p$). That is $f, g : M \to M$ are considered equal if $f(x)
= g(x) \in \mathfrak p$. The nonlogical symbols of $L$ are interpreted as
follows. Given an $n$-ary relation $R$ and elements $f_1, ..., f_n$ of $M^*$,
the formula $R(f_1, ..., f_n)$ holds in $M^*$ if and only if $R(f_1(x), ...,
f_n(x)) \in \mathfrak p$. Just as in the regular ultrapower construction, one
sees that this does not depend on the representatives $f_1, ..., f_n$. Similar
definitions are applied to the interpretation of constant and functional
symbols. We have the analogue of \L{}o\'{s}'s Theorem.

\begin {lemma}
Given an $L$-formula $\phi(y_1, ..., y_n)$ and elements $f_1, ..., f_n \in M^*$
we have $M^* \models \phi(f_1, ..., f_n)$ if and only if $\phi(f_1(x), ...,
f_n(x)) \in \mathfrak p$.
\end {lemma}

\begin {proof}
By induction of $\phi$. The only nontrivial case is when $\phi$ is of the form
$\exists y_0 \psi(y_0, y_1, ..., y_n)$. Assuming $M^* \models \exists y_0
\psi(y_0, f_1, ..., f_n)$, there exists $f_0 \in M^*$ such that $M^* \models
\psi(f_0, f_1, ..., f_n)$. Then by the induction hypothesis $\psi(f_0(x),
f_1(x), ..., f_n(x)) \in \mathfrak p$. And since it implies $\exists y_0
\psi(y_0, f_1(x), ..., f_n(x))$, the latter is also in $\mathfrak p$. Conversely
assume that $\exists y_0 \psi(y_0, f_1(x), ..., f_n(x)) \in \mathfrak p$. Since
$T$ has definable Skolem functions, there is a definable function $f_\psi : M
\to M$ such that $\forall x (\exists y_0 \psi(y_0, f_1(x), ..., f_n(x)) \to
\psi(f_\psi(x), f_1(x), ..., f_n(x)))$. Therefore $\psi(f_\psi(x), f_1(x), ...,
f_n(x)) \in \mathfrak p$ and by the induction hypothesis $M^* \models \exists
y_0 \psi(y_0, f_1, ..., f_n)$.
\end {proof}

\begin {corollary}
If we identify each element of $M$ with the constant function in $M^*$, then
$M^*$ is an elementary extension of $M$.
\end {corollary}

Using the fact that $\mathfrak p$ is a strongly regular type, the result of the
previous Lemma can also be expressed as follows. Assume that $a$ is generic over
parameters defining $f_1, ..., f_n$ in $M$, then $M^* \models \phi(f_1, ...,
f_n)$ if and only if $M \models \phi(f_1(a), ..., f_n(a))$.

We call $M^*$ the canonical extension of $M$. (This explains the terminology
self-extending.) The following property of the canonical extension is crucial
for our purposes.

\begin {proposition}
If $M$ is an infinite dimensional model, then its canonical extension $M^*$ is a
proper extension where every new element realises $\mathfrak p|_M$.
\end {proposition}

\begin {proof}
First note that since $\mathfrak p|_M$ is not isolated, the identity function in
$M^*$ is different from all constant functions modulo $\mathfrak p$. Hence $M^*$
is a proper extension of $M$.

On the other hand, let $f \in M^*$ be an element that does not realise
$\mathfrak p|_M$. Then there is an $L(M)$-formula $\psi(\bar a, y) \not \in
\mathfrak p$ such that $M^* \models \psi(\bar a, f)$. Let $\phi(\bar b, x, y)$
define $f$ in $M$. Pick an element $c$ generic over $\bar a \bar b$. By the
assumption $d = f(c)$ satisfies $\psi(\bar a, y)$ in $M$. But then $d \in
\cl_{\mathfrak p}(\bar a)$ and since $c$ is generic over $\bar a \bar b$ we see
that $f(x) = d$ for every generic $x$. Thus $f$ is a constant function (modulo
$\mathfrak p$) and so is already in $M$.
\end {proof}

\begin {theorem}
Let $\mathfrak p$ be a symmetric strongly regular type in a theory $T$ with
definable Skolem functions. Then $T$ has quasiminimal models of arbitrarily
large cardinalities.
\end {theorem}

\begin {proof}
We follow the approach outlined in the beginning of the section. Let $\kappa$ be
an uncountable cardinal and let $A = (a_\alpha : \alpha < \kappa)$ be a Morley
sequence in $\mathfrak p$. We show that $\dcl(A)$ is quasiminimal. Let $A_0
\subseteq A$ be countably infinite. Let $M$ be a countable model containing
$A_0$. Iterating the previous Proposition we obtain an extension $N$ such that
$\dim_{\mathfrak p}(N) = \kappa$ and $\cl_{\mathfrak p}(A_0) \cap N =
\cl_{\mathfrak p}(A_0) \cap M$ is countable. By conjugating $N$ with an
automorphism if necessary we may assume that $A \subseteq N$. But then $\dcl(A)
\subseteq N$ and hence $\cl_{\mathfrak p} \cap \dcl(A)$ is also countable. This
shows that $\dcl(A)$ is quasiminimal.
\end {proof}

In the rest of this section we do not assume that $T$ has definable Skolem
function. Instead we add functional symbols to the language that will resemble
Skolem functions. We use this to prove that if $T$ has a quasiminimal model of
cardinality $\beth_{\omega_1}$, then it has quasiminimal models of all
uncountable cardinalities. The existence of such a cardinal follows from a very
general argument due to Hanf, and so the least such cardinal is often called the
Hanf number. In this terminology we show that the Hanf number of the existence
of quasiminimal models is at most $\beth_{\omega_1}$.

\begin {theorem}
If $T$ has a quasiminimal model of cardinality $\beth_{\omega_1}$, then it has
quasiminimal models of arbitrarily large cardinalities.
\end {theorem}

\begin {proof}
Let $M \models T$ with $|M| = \beth_{\omega_1}$. Then by the Fact \ref
{pilovic-quasi} (after adding parameters) there is a global symmetric regular
type $\mathfrak p$ that extends the type of cocountable subsets of $|M|$. Since
$M$ is quasiminimal, its dimension (with respect to $\cl_{\mathfrak p}$) is
$\beth_{\omega_1}$ and so there is a Morley sequence $(a_\alpha : \alpha <
\beth_{\omega_1})$ in $\mathfrak p$.

Expand the language by $n$-arey functional symbols $f_n^i$ for each $i, n <
\omega$ (nullary functional symbols being constant symbols). Interpret $f_n^i$
in $M$ such that for every $\alpha_1 < ... <\alpha_n < \kappa$ the set
$\{f_n^i(a_{\alpha_1}, ..., a_{\alpha_n}) : i < \omega\}$ enumerates the closure
$\cl_{\mathfrak p}(a_{\alpha_1}, ...., a_{\alpha_n})$ in $M$. Denote the
resulting language, theory and structure by $L'$, $T'$ and $M'$ respectively.

Given a cardinal $\kappa$, there is an indiscernible sequence $B = (b_\beta :
\beta < \kappa)$ in the monster model $\mathfrak C'$ of $T'$ such that for every
$m$ there are $\alpha_1 < ... < \alpha_m$ satisfying
$$\tp^{\mathfrak C'}(b_0,...,b_{m-1}) =
\tp^{M'}(a_{\alpha_1},...,a_{\alpha_m}).$$
(This is commonly known as ``Morley's method''. A proof can be found for example
in \cite {tent-ziegler} or \cite {casanovas} where the length of the sequence is
assumed to be $\beth_{(2^{|T|})^+}$. The better bound of $\beth_{\omega_1}$ for
countable theories is from \cite {simple-primer}.)

Now let $N = \{f_n^i(b_{\beta_1}, ..., b_{\beta_n}) : i, n < \omega, \beta_1
< ... < \beta_n < \kappa\} \subseteq \mathfrak C'$ be the closure of this
indiscernible sequence under all $f_n^i$-s. Note that $B \subset N$ is a Morley
sequence in $\mathfrak p$.  Also each $f_n^i(b_{\beta_1}, ..., b_{\beta_n}) \in
\cl_{\mathfrak p}(b_{\beta_1}, ..., b_{\beta_n})$ (since this is encoded in the
type of the corresponding sequence in $M'$). Hence $N \subseteq \cl_{\mathfrak
p}(B)$ has dimension $\kappa$.

We claim that $N$ is a quasiminimal model of $T$. This essentially follows from
the fact that for every $\beta_1 < ... < \beta_n < \kappa$ there are $\alpha_1 <
... < \alpha_n < \beth_{\omega_1}$ such that the mapping sending
$f_n^i(b_{\beta_1},...,b_{\beta_n}) \mapsto
f_n^i(a_{\alpha_1},...,a_{\alpha_n})$ is elementary in $L$. Now by embedding a
suitable fragment of $N$ inside $M$ in this way we see that the elements of $N
\setminus \{f_n^i(b_{\beta_1}, ..., b_{\beta_n}) : i < \omega\}$ are generic
over $b_{\beta_1}, ..., b_{\beta_n}$.  This shows that $N$ is quasiminimal.
Similarly each satisfiable $L$-formula over $f_n^{i_1}(b_{\beta_1}, ...,
b_{\beta_n}), ..., f_n^{i_m}(b_{\beta_1}, ..., b_{\beta_n})$ is either in
$\mathfrak p$ and hence realised in $B$ (since $B$ is infinite dimensional), or
not in $\mathfrak p$ and realised in $\{f_n^i(b_{\beta_1}, ..., b_{\beta_n}) : i
< \omega\}$. This shows that $N$ is an $L$-elementary substructure of $\mathfrak
C'$ and hence a model of $T$.
\end {proof}

\section {Stable Theories}

In this section we assume that apart from having a strongly regular type
$\mathfrak p$, the theory $T$ is stable. This automatically makes $\mathfrak p$
symmetric (all indiscernible sequence are totally indiscernible in stable
theories). And since all symmetric global regular types are
$\emptyset$-definable, the results of Section \ref {gensec} apply here.

Some standard background on stable theories would be assumed. To construct a
quasiminimal model we use the technology of local isolation and local atomicity.
It was first introduced by \cite {lachlan} in order to prove a two cardinal
theorem for stable theories. Here we present the background on local isolation
in order to make our account more complete. More details can be found in
advanced books on stability theory such as \cite {shelah, baldwin-stab}. Our
notion of local isolation coincides with $\mathbf F^l_{\aleph_0}$-isolation of
\cite {shelah}.

\begin {definition}
A type $p(\bar x) \in S(A)$ is called {\em locally isolated} if for every
$L$-formula $\phi(\bar x, \bar y)$ there is a formula $\psi_\phi(\bar x) \in p$
such that $\psi_\phi(\bar x) \vdash p(\bar x)|_\phi$. (Here $p(\bar x)|_\phi =
\{\phi(\bar x, \bar a) : \phi(\bar x, \bar a) \in p\} \cup \{\lnot \phi(\bar x,
\bar a) : \lnot \phi(\bar x, \bar a) \in p\}$ is the $\phi$-part of $p$.) A set
$B$ is called {\em locally atomic} (over $A$) if every type (over $A$) realised
in $B$ is locally isolated.
\end {definition}

The notion of local isolation can be seen through the topology of $\phi$-types
as follows. Let $S_\phi(A)$ be the space of all complete $\phi$-types over $A$.
Topologise $S_\phi(A)$ by taking the basic clopen sets to be of the form $[\psi]
= \{p \in S_\phi(A) : p \vdash \psi\}$ where $\psi$ is any boolean combination
of $\phi$-formulas over $A$. Then $S_\phi(A)$ is a boolean topological space
(compact, Hausdorff and totally disconnected) and the canonical restriction
$\sigma_\phi : S(A) \to S_\phi(A)$ is continuous (and hence also closed).  Now a
type $p \in S(A)$ will be locally isolated if and only if
$\sigma_\phi^{-1}(\sigma_\phi(p))$ is a neighbourhood of $p$ for every $\phi$
(i.e. contains an open set containing $p$).

Recall that in stable theories $S_\phi(\mathfrak C)$ is a scattered topological
space of finite Cantor-Bendixson rank (see \cite {casanovas}). For a partial
type $p$ we denote by $\CB_\phi(p)$ and $\Mlt_\phi(p)$ the Cantor-Bendixson rank
and degree of the set $\sigma_\phi([p]) = \{q \in S_\phi(\mathfrak C) : p \cup q
\text { is consistent}\}$ in $S_\phi(\mathfrak C)$.

We would like to show that for every subset $A$, there is a model $M \supseteq
A$ that is locally atomic over $A$. For that we prove

\begin {lemma}
For any set $A$, locally isolated types are dense in $S(A)$.
\end {lemma}

\begin {proof}
Given a formula $\chi(\bar x)$ over $A$, we need to find a locally isolated type
that contains $\chi$. Enumerate $(\phi_n(\bar x, \bar y_n) : n < \omega)$ all
formulas in $L$. Put $p_0 = \{\chi(\bar x)\}$ and $p_{n+1} = p_n \cup
\{\psi_n(\bar x, \bar a_n)\}$ where $\psi_n$ is chosen such that
$\CB_{\phi_n}(p_{n+1})$ is the least possible and among those
$\Mlt_{\phi_n}(p_{n+1})$ is the least possible. We claim that $p = \bigcup_{n <
\omega} p_n$ has a unique completion which is locally isolated.  Indeed given
$\phi_n(\bar x, \bar y_n)$ and $\bar b_n \in A$, we can't have both $\phi_n(\bar
x, \bar b_n)$ and $\lnot \phi_n(\bar x, \bar b_n)$ consistent with $p_{n+1}$ as
it contradicts the choice of $\psi_n$. Hence $p_{n+1}$, which is finite isolates
the $\phi_n$-part of the completion of $p$.
\end {proof}

This allows us to iteratively realise formulas by locally isolated types similar
to the constructible models in $\omega$-stable theories. To make the whole
construction locally atomic over $A$ we need

\begin {lemma}
\label {local-trans}
If $A$ is locally atomic over $BC$ and $B$ is locally atomic over $C$, then $AB$
is locally atomic over $C$.
\end {lemma}

\begin {proof}
Let $\bar a \in A$ and $\bar b \in B$ and $\phi(\bar x, \bar y, \bar z)$ be a
formula. It is enough to find a formula $\psi(\bar x, \bar y) \in \tp(\bar a\bar
b/C)$ such that for every $\bar c \in C$ for which $\phi(\bar a, \bar b, \bar
c)$ holds we have $\psi(\bar x,\bar y) \to \phi(\bar x, \bar y, \bar c)$. Since
$\tp(\bar a/BC)$ is locally isolated, we have a formula $\chi(\bar x, \bar b')$
over $C$ in $\tp(\bar a/BC)$ such that $\chi(\bar x, \bar b') \to \phi(\bar x,
\bar b, \bar c)$ whenever $\phi(\bar a, \bar b, \bar c)$ holds. Now we have
$\forall \bar x (\chi(\bar x, \bar z) \to \phi(\bar x, \bar y, \bar c)) \in
\tp(\bar b\bar b'/C)$. By the local isolation of the latter there is a formula
$\sigma(\bar y, \bar z) \in \tp(\bar b'\bar b/C)$ such that
$$\sigma(\bar y, \bar z) \to \forall \bar x(\chi(\bar x, \bar z) \to \phi(\bar
x, \bar y, \bar c))$$ 
whenever $\forall x(\chi(\bar x, \bar b') \to \phi(\bar x, \bar b, \bar c))$
holds. Or equivalently
$$(\exists \bar z (\sigma(\bar y, \bar z) \land \chi(\bar x, \bar z))) \to
\phi(\bar x, \bar y, \bar c).$$
Thus $\exists \bar z(\sigma(\bar y, \bar z) \land \chi(\bar x, \bar z))$ is the
required formula $\psi(\bar x, \bar y)$.
\end {proof}

Now the two combine to give

\begin {proposition}
For every $A$ in a stable theory there is a model $M \supseteq A$ that is
locally atomic over $A$ and $|M| = |A| + \aleph_0$.
\end {proposition}

\begin {proof}
We imitate the construction of prime models in $\omega$-stable theories. Given
$B$ construct $B'$ as follows: enumerate $(\phi_\alpha(\bar x_\alpha) : \alpha <
\kappa)$ all consistent $L(B)$ formulas that are not realised in $B$. Then add
realisations $b_\alpha$ to $B'$ such that $\tp(b_{\alpha}/B\{b_\beta : \beta <
\alpha\})$ is locally isolated. By Lemma \ref {local-trans}, $B'$ is locally
atomic over $B$. Now take $A_0 = A$ and $A_{n+1} = A_n'$. Finally $M =
\bigcup_{n < \omega} A_n$ is the required model.
\end {proof}

The final tool that we need for working with local isolation is the analogue of
the Open Mapping Theorem.

\begin {proposition} (Local Open Mapping Theorem)
If $p \in S(A)$ does not fork over $B \subseteq A$ and is locally isolated, then
so is $p|_B$.
\end {proposition}

\begin {proof}
Let $N(A/B)$ be the set of types in $S(A)$ that don't fork over $B$. By the
usual Open Mapping Theorem the restriction map $\pi : N(A/B) \to S(B)$ is open.
Now fix a formula $\phi(\bar x, \bar y)$. By the local isolation of $p$ we have
that $\sigma_\phi^{-1}(\sigma_\phi(p)) \cap N(A/B)$ is a neighbourhood of $p$ in
$N(A/B)$. Hence its image under $\pi$, which is contained in
$\sigma_\phi^{-1}(\sigma_\phi(p|_B))$ is a neighbourhood of $p|_B$.
\end {proof}

We can use locally atomic models to construct a proper extension of a given
model $M$ where all the new elements realise the type $\mathfrak p|_M$. (In the
previous section we used the canonical extension of a model for theories with
definable Skolem functions for this purpose.) Indeed let $a$ realise $\mathfrak
p|_M$ and let $N$ be locally atomic over $Na$. Now let $b \in N \setminus M$.
Then $\tp(b/Ma)$ is locally isolated and $\tp(b/M)$ is not (otherwise it would
be realised). Hence by the Local Open Mapping Theorem $b \forks_M a$. But then
by symmetry we have $a \forks_M b$, which means that $a \not \models \mathfrak
p|_{Mb}$, i.e. $a \in \cl_{\mathfrak p}(Mb)$. Hence $b \not \in \cl_{\mathfrak
p}(M)$, as otherwise $a \in \cl_{\mathfrak p}(M)$.

However the situation is a bit more delicate here. The problem is that there is
no equivalent notion of local primness. That is we cannot in general embed a
locally atomic or a locally constructible model over $A$ in an arbitrary model
containing $A$. So we cannot use the same idea as in the previous section to
construct arbitrarily large quasiminimal models. 

As an alternative we could try to construct a locally atomic model over a Morley
sequence $A$ in the following way. Enumerate $A = (a_\alpha : \alpha < \kappa)$.
Then build a sequence $(M_\alpha : \alpha < \kappa)$ such that $M_{\alpha+1}$ is
locally atomic over $M_\alpha a_\alpha$.  Then all the elements in $M_{\alpha+1}
\setminus M_\alpha$ will be generic over $M_{\alpha}$, so we do not extend the
closure. The problem with this, however, is that for $\alpha \ge \aleph_1$ the
models $M_\alpha$ will be uncountable. And so in general $M_{\alpha+1}$ adds
uncountably many elements.

To remedy this we employ a construction along the lines of Shelah's
Excellence/NOTOP. The technique was adapted for local atomicity in \cite
{bays-covers}. Some of our arguments are adapted from that paper. First let us
introduce some notation.

\begin {definition}
Let $I$ be a downward-closed set of subsets (i.e. $s \subseteq t \in I$ implies
$s \in I$). An {\em $I$-system} is a collection $(M_s : s \in I)$ of models such
that $M_s \preceq M_t$ whenever $s \subseteq t$. If $J \subseteq I$ we denote
$M_J = \bigcup_{s \in J} M_s$.

For $s \in I$ denote $<s := \mathcal P(s) \setminus \{s\}$ and $\not \geq s :=
\{t \in I : t \not \supseteq s\}$.

The system is {\em independent} if $M_s \nforks_{M_{<s}} M_{\not \geq s}$.

An {\em enumeration} of $I$ is an ordering $(s_\alpha : \alpha \in \kappa)$ of
$I$ such that $\alpha \leq \beta$ whenever $s_\alpha \subseteq s_\beta$. If an
enumeration is fixed we write $< \alpha$ for $\{s_\beta: \beta < \alpha\}$. I.e.
$M_{< \alpha} = \bigcup_{\beta < \alpha} M_{s_\beta}$.
\end {definition}

\begin {definition}
Given two subsets $A \subseteq B$ we say that $A$ is Tarski-Vaught in $B$ (in
symbols $A \subseteq_{TV} B$) if every formula over $A$ realised in $B$ is
already realised in $A$.
\end {definition}

We use the Tarski-Vaught condition to lift local isolation to larger sets.

\begin {proposition}
If $A \subseteq_{TV} B$ and $\bar c$ is locally isolated over $A$, then it is
locally isolated over $B$.
\end {proposition}

\begin {proof}
Let $\phi(\bar x, \bar y)$ be a formula without parameters and assume that
$\psi(\bar x)$ isolates the $\phi$-part of $\tp(\bar c/A)$. Then $\psi(\bar x)$
isolates a $\phi$-type over $B$. Indeed if for some $\bar b \in B$ we have
$\exists \bar x (\phi(\bar x) \land \phi(\bar x, \bar b)) \land \exists \bar x
(\phi(\bar x) \land \lnot \phi(\bar x, \bar b))$, then the same should be true
for some $\bar a \in A$. Thus $\psi(\bar x)$ isolates the $\phi$-part of
$\tp(\bar c/B)$.
\end {proof}

Finally we need the following Lemma from \cite {shelah}.

\begin {lemma} [{TV Lemma, \cite [XII.2.3(2)] {shelah}}]
Let $(M_s : s \in I)$ be an independent system in a stable theory. Assume that
$J \subseteq I$ is such that for all $s \in I$ if $s \subseteq \bigcup J$, then
$s \in J$. Then $M_J \subseteq_{TV} M_I$.
\end {lemma}

\begin {theorem}
\label {stable-main}
Let $\mathfrak p$ be a strongly regular type in a stable theory $T$. Then $T$
has quasiminimal models of arbitrarily large cardinalities.
\end {theorem}

\begin {proof}
Let $I$ be the set of all finite subsets of some uncountable cardinal $\kappa$.
We inductively build an independent $I$-system of countable models as follows.
Enumerate $I = \{s_\alpha : \alpha < \kappa\}$ so that $s_\alpha \subseteq
s_\beta$ implies $\alpha \leq \beta$. Note that this implies $s_0 = \emptyset$.
For $M_\emptyset$ pick a countable infinite dimensional model. Given an ordinal
$\alpha$, assume that $M_{s_\beta}$ has been constructed for $\beta < \alpha$.
If $s_\alpha = \{\delta\}$ is a singleton, pick $a_\delta$ a realisation of
$\mathfrak p|_{M_{< \alpha}}$ and let $M_{\{\delta\}}$ be a countable locally
atomic model over $M_\emptyset a_\delta$ that is independent from $M_{< \alpha}$
over $M_\emptyset$. Otherwise let $M_s$ be a countable locally atomic model over
$M_{<s}$ that is independent from $M_{< \alpha}$ over $M_{<s}$. Then $(M_s : s
\in I)$ is an independent system (this is \cite [Lemma XII.2.3(1)] {shelah},
whose proof is ``an exercise in non-forking'').

Now $M_I$ is a model (by Tarski-Vaught test) and $|M_I| = \kappa$. We claim that
$M_I$ is quasiminimal. Fix $s = \{\alpha_1, ..., \alpha_n\} \subset \kappa$. We
claim that all elements of $M_I \setminus M_s$ realise $\mathfrak p |_{M_s}$
over $M_s$. Since the later is countable, this implies that every subset of
$M_I$ definable over $M_s$ is either countable or cocountable. Since $s$ was
arbitrary we conclude that $M_I$ is quasiminimal.

Let $A = \{a_\alpha : \alpha < \kappa\}$. We show that $M_I$ is locally atomic
over $M_sA$. With the above enumeration of $I$ we show by induction on $\alpha$
that $M_{< \alpha}$ is locally atomic over $M_sA$. This is clear if $\alpha = 0$
or $\alpha$ is a limit ordinal. So assume that it holds for $\alpha$ and let us
prove it for $\alpha+1$. Consider several cases.
\begin {itemize}
\item If $s_\alpha \subseteq s$, then clearly $M_{s_\alpha}$ is locally isolated
	over $M_{< \alpha}M_sA$.
\item If $|s_\alpha| > 1$, then $M_{s_\alpha}$ is locally atomic over $M_{<
	s_\alpha}$. By the first clause we may assume that $s_\alpha \not
	\subseteq s$. But by TV Lemma (with $I = \{s_\beta : \beta < \alpha\}
	\cup \mathcal P(s)$ and $J = \mathcal P(s_\alpha) \setminus
	\{s_\alpha\}$) we have $M_{< s_\alpha} \subseteq_{TV} M_{
	< \alpha}M_s$. Also $M_{< \alpha}M_s \subseteq_{TV} M_{< \alpha}M_sA$
	since $M_\emptyset$ is infinite dimensional. Hence $M_{s_\alpha}$ is
	locally atomic over $M_{< \alpha}M_sA$.
\item If $s_\alpha = \{\delta\}$ is a singleton, then $M_{s_\alpha}$ is locally
	atomic over $M_\emptyset a_\delta$. By the first clause we may assume
	that $\delta \not \in s$ and hence $a_\delta$ realises $\mathfrak p$
	over $M_{< \alpha}M_s$. Now since $M_\emptyset$ is a model, we have
	$M_\emptyset \subseteq_{TV} M_{< \alpha}M_s$. Using the definability of
	$\mathfrak p$ we can show that $M_\emptyset a_\delta \subseteq_{TV} M_{<
	\alpha}M_sa_{\delta}$ and as before $M_{< \alpha}M_sa_{\delta}
	\subseteq_{TV} M_{< \alpha}M_sA$. Thus we conclude that $M_{s_\alpha}$
	is locally atomic over $M_{< \alpha}M_sA$.
\end {itemize}
Thus in all cases we have $M_{s_\alpha}$ is locally atomic over $M_{<
\alpha}M_sA$.  But by the induction hypothesis $M_{< \alpha}$ is locally atomic
over $M_sA$. Hence $M_{s_\alpha}M_{< \alpha} = M_{< \alpha+1}$ is locally atomic
over $M_sA$. This completes the induction.

Now given $b \in M \setminus M_s$, it is locally atomic over $M_sA$. Since $b$
is not locally atomic over $M_s$ (otherwise its type would be realised in
$M_s$), we conclude by the Local Open Mapping Theorem that $b \not \nforks_{M_s}
A$.  Pick a finite subset $A' \subset A$ such that $b \not \nforks_{M_s} A'$. We
can also assume that $A'$ is disjoint from $M_s$. Now by symmetry we have $A'
\not \nforks_{M_s} b$. Since $A'$ is a Morley sequence in $\mathfrak p$ over
$M_s$, we conclude that $\dim(A'/M_sb) < \dim(A'/M_s)$ (where $\dim$ is in the
sense of pregeometry $\cl_{\mathfrak p}$). Hence $b \not \in \cl_{\mathfrak
p}(M_s)$.  This finishes the proof.
\end {proof}

\section {$\omega$-stable Theories}

In this section we assume that the theory $T$ is $\omega$-stable. This allows us
to construct prime models over every subset. It is easy to see that a prime
model over an uncountable Morley sequence in $\mathfrak p$ must be quasiminimal.
Here we show more: the class $\mathcal C$ of prime models over Morley sequences
in $\mathfrak p$ is a quasiminimal excellent class (see Definition \ref
{qpg-class}). Since $\mathcal C$ is clearly uncountably categorical, its
excellence was expected and possibly known to experts. But we are not aware of a
published proof. 

Quasiminimal excellent classes play an important role in nonelementary
categoricity. They were originally introduced in \cite {zilber-qme} where it is
proven that a quasiminimal excellent class is uncountably categorical. The
original formulation contains a technical axiom called {\em excellence}. It was
thought to be the key to categoricity, until \cite {bhhkk} showed that
excellence follows from the rest of the axioms (see also a direct proof of
categoricity in \cite {myself}). \cite {bhhkk} called an infinite dimensional
structure in such a class a {\em quasiminimal pregeometry structure} and
following this we called the entire class {\em quasiminimal pregeometry class}
in \cite {myself}.  However in this paper we are dealing with structure that are
quasiminimal and have a pregeometry but are not quasiminimal pregeometry
structures in the sense of \cite {bhhkk}. So we have reverted the terminology
back to quasiminimal excellent. 

The following simple observation will be used repeatedly, so it is worth stating
explicitly. We have already used a variant of it for locally atomic model in the
previous section.

\begin {proposition}
\label {prime-ind}
Let $M$ be a model $a \not \in \cl_{\mathfrak p}(M)$ and $N$ prime over $Ma$.
Then for every $b \in N \setminus M$ we have $a \in \cl_{\mathfrak p}(Mb)$. 
Hence by the exchange property also $b \in \cl_{\mathfrak p}(Ma) \setminus
\cl_{\mathfrak p}(M)$.
\end {proposition}

\begin {proof}
Indeed since $\tp(b/Ma)$ is isolated, whereas $\tp(b/Ma)$ is not, by the Open
Mapping Theorem we have $b \not \nforks_M a$. Hence $a \not \nforks_M b$. So
that $a \in \cl_{\mathfrak p}(Mb)$.
\end {proof}

Let us first show that each uncountable model in $\mathcal C$ is quasiminimal.
One consequence of quasiminimality is that there are no uncountable
indiscernible sequences except in $\mathfrak p$. This will help us establish
primeness of models in some cases.

\begin {lemma}
If $A$ is an uncountable Morley sequence in $\mathfrak p$, then the prime model
over $A$ is quasiminimal.
\end {lemma}

\begin {proof}
This follows the standard pattern we have used already in previous sections.
For a fixed countable $A_0 \subseteq A$ we can construct a model $N$ containing
$A$ such that $\cl_{\mathfrak p}(A_0) \cap N$ is countable. To do so enumerate
$A \setminus A_0 = (a_\alpha : \alpha < \kappa)$ and take $N_{\alpha+1}$ to be
prime over $N_\alpha a_\alpha$. Let $N_0$ be an arbitrary countable model
containing $A_0$ and $N_\delta = \bigcup_{\alpha < \delta} N_\delta$ for a limit
ordinal $\delta$.  Then take $N = N_\kappa$. Finally $M$ embeds into $N$ over
$A$ showing that $\cl_{\mathfrak p}(A_0) \cap M$ is also countable.
\end {proof}

In order to satisfy some technical conditions of Definition \ref {qpg-class}, we
make two further assumptions, which can be achieved by expanding the language.
Firstly we assume that the theory $T$ has quantifier elimination. Otherwise we
can consider its Morleysation, i.e. its expansion with predicate symbols for all
$\emptyset$-definable relations. Secondly we assume that $\mathfrak p|_A$ is not
isolated for all finite $A$. (Note that this implies that $\mathfrak p|_A$ is
not isolated for all $A$.) Otherwise we can add a countably infinite Morley
sequence in $\mathfrak p$ to the language. The result of the second assumption
is

\begin {lemma}
\label {fin-iso}
\begin {enumerate}
\item If $M$ is a model and $A \subseteq M$, then $\cl_{\mathfrak p}(A) \cap M$
	is a model.
\item If $M$ is prime over a Morley sequence $A$ in $\mathfrak p$, then $A$ is
	a basis for $M$ (i.e. $M \subseteq \cl_{\mathfrak p}(A)$).
\end {enumerate}
\end {lemma}

\begin {proof}
\begin {enumerate}
\item Let $N = \cl_{\mathfrak p}(A) \cap M$. We use the Tarski-Vaught test. Let
	$\phi(x, \bar b)$ be a formula over $N$. Since $\mathfrak p|_{\bar b}$
	is not isolated, there is a consistent formula $\psi(x, \bar b)$ that
	implies $\phi$ and is not in $ \mathfrak p$. Now every element in $M
	\setminus N$ realises $\mathfrak p|_{\bar b}$. Hence a realisation of
	$\psi$ must be in $N$.
\item Given $b \in M$, we have that $\tp(b/A)$ is isolated. Hence $b \not
	\models \mathfrak p|_A$ and so $b \in \cl_{\mathfrak p}(A)$.
\end {enumerate}
\end {proof}

Note that both assertions can fail without assuming that $\mathfrak p|_A$ is
isolated for all finite $A$. For example both assertions fail in the theory of
an infinite set in the empty language, where $\mathfrak p$ is the unique
non-algebraic type (so the theory is strongly minimal).

Now get prepared for a large definition. In the following a partial embedding is
a partial map that preserves quantifier free formulas. 

\begin {definition}
\label {qpg-class}
A {\em quasiminimal excellent} class is a collection $\mathcal C$ of pairs $(H,
\cl_H)$ where $H$ is a structure and $\cl_H$ is a pregeometry on $H$ satisfying
the following conditions.
\begin {enumerate}
\item {\em Closure under isomorphisms}\newline
	If $(H, \cl_H) \in \mathcal C$ and $f : H \to H'$ is an
	isomorphism, then $(H', \cl_{H'}) \in \mathcal C$, where
	$\cl_{H'}$ is defined as $\cl_{H'}(X') = f(\cl_H(f^{-1}(X')))$ for $X'
	\subseteq H'$.

\item {\em Quantifier free theory}\newline
	If $(H, \cl_H), (H', \cl_H') \in \mathcal
	C$, then $H$ and $H'$ satisfy the same quantifier free sentences. 

\item {\em Pregeometry}
\begin {itemize}
\item For each $(H, \cl_H) \in \mathcal C$ the closure of any
	finite set is countable.
\item If $(H, \cl_H) \in \mathcal C$ and $X \subseteq H$, then
	$\cl_H(X)$ is a substructure of $H$ and together with the restriction of
	$\cl_H$ it is in $\mathcal C$.
\item If $(H, \cl_H), (H', \cl_{H'}) \in \mathcal
	C$, $X \subseteq H$, $y \in H$ and $f : H \to H'$ is a partial
	embedding defined on $X \cup \{y\}$, then $y \in \cl_H(X)$ if and only
	if $f(y) \in \cl_{H'}(f(X))$.
\end {itemize}

\item {\em Uniqueness of the generic type over countable closed models}\newline
	Let $(H, \cl_H), (H', \cl_{H'}) \in \mathcal
	C$, subsets $G \subseteq H, G' \subseteq H'$ be countable closed or
	empty and $g : G \to G'$ be an isomorphism.  If $x \in H, x' \in H'$ are
	independent from $G$ and $G'$ respectively, then $g \cup \{(x,
	x')\}$ is a partial embedding.

\item {\em $\aleph_0$-homogeneity over countable closed models}\newline
	Let $(H, \cl_H), (H', \cl_{H'}) \in \mathcal
	C$, subsets $G \subseteq H, G' \subseteq H'$ be countable closed or
	empty and $g : G \to G'$ be an isomorphism.  If $g \cup f : H \to H'$
	is a partial embedding, $X = \dom(f)$ is finite and $y \in \cl_H(X \cup
	G)$, then there is $y' \in H'$ such that $g \cup f \cup \{(y,
	y')\}$ is a partial embedding.
\end {enumerate}
\end {definition}

We don't elaborate on the definition any further. The interested reader can
consult \cite {zilber-qme}, \cite {bhhkk} or \cite {myself}. Note however that
given a strongly regular type $\mathfrak p$ in an arbitrary theory, the class of
elementary submodels of the monster model $\mathfrak C$ (closed under
isomorphisms) together with the restriction of $\cl_{\mathfrak p}$ satisfies
axioms 1 and 2. If further the theory has quantifier elimination and $\mathfrak
p|_A$ is not isolated for all finite $A$, then the class satisfies axioms 3 and
4, except the countable closure property. To satisfy the countable closure
property we need to take the class of elementary submodels of a large
quasiminimal structure (if such a structure exists). Finally satisfying axiom 5
is the real challenge. For that we need to show that in a prime model over a
Morley sequence, if we pick a different Morley sequence, then the respective
closure is prime over it. The following lemma is extracted from \cite {makkai}.

\begin {lemma}
\label {makkai}
Assume $M$ is a model $a \models \mathfrak p|_M$ and $N$ is a prime model over
$Ma$. Let $b \in N \setminus M$. Then $N$ is prime over $Mb$. Consequently there
is an automorphism in $\Aut(N/M)$ taking $a$ to $b$.
\end {lemma}

\begin {proof}
By Proposition \ref {prime-ind} we have $a \not \nforks_M b$. Let $\theta(x, b)
\in \tp(a/Mb)$ be a formula over $Mb$ such that any type containing it forks
over $M$. Let $\phi(a, y)$ be a formula over $Ma$ that isolates $\tp(b/Ma)$. We
claim that $\theta(x, b) \land \phi(x, b)$ isolates $\tp(a/Mb)$.  It is enough
to show that for every $a' \in N$ such that $\theta(a', b) \land \phi(a', b)$ we
have $a \equiv_{Mb} a'$. Given such an $a'$ we have $a' \not \nforks_M b$. Then
$a' \in N \setminus M$. Hence by Proposition \ref {prime-ind} again $a' \models
\mathfrak p_M$. Now we can get the desired conclusion from $\phi(a', b)$. Indeed
let $f \in \Aut(\mathfrak C/M)$ be an automorphism that maps $a'$ to $a$. Then
$\phi(a, f(b))$ holds. Since $\phi(a, y)$ isolates $\tp(b/Ma)$ we have $ba
\equiv_M f(b)a \equiv_M ba'$. Hence $a \equiv_{Mb} a'$.

Now since $N$ is constructible over $Ma$ and $\tp(a/Mb)$ is isolated, we
conclude that $N$ is constructible and hence prime and atomic over $Mb$.
\end {proof}

\begin {proposition}
\label {qpg-prime}
Assume that $M$ is prime over a Morley sequence $A = (a_\alpha : \alpha <
\kappa)$ in $\mathfrak p$. Let $B = (b_\alpha : \alpha < \mu) \subseteq M$ is
another Morley sequence in $\mathfrak p$. Then $\cl_{\mathfrak p}(B) \cap M$ is
prime over $B$.
\end {proposition}

\begin {proof}
By Lemma \ref {fin-iso}, the closure $\cl_{\mathfrak p}(B) \cap M$ is a model.
Also by the countable closure property it does not contain an uncountable
indiscernible sequence over $B$. Thus it remains to show that $\cl_{\mathfrak
p}(B) \cap M$ is atomic over $B$.

We can assume that $B$ is finite. Indeed every $a \in \cl_{\mathfrak p}(B) \cap
M$ is in the closure of a finite subset $B_0 \subseteq B$. Assume we prove that
$a$ is atomic over $B_0$. Then given $b \in B \setminus B_0$, we have that
$\tp(a/B_0) \vdash \tp(a/B_0b)$ (see the discussion after Definition \ref
{sreg-type}). Hence $\tp(a/B_0b)$ is also isolated. Iterating this we obtain
that $\tp(a/B)$ is isolated. So assume that $B = \{b_0, ..., b_{n-1}\}$ is
finite.

Further we may assume that $B$ is a subset of $A$. Indeed by the exchange
property there is $a \in A$ such that for $A' = A \setminus \{a\}$ we have that
$A' \cup \{b_0\}$ is a basis of $M$. Now $N = \cl_{\mathfrak p}(A') \cap M$ is a
model by Lemma \ref {fin-iso}. Further $M$ is prime over $Na$ (it is atomic over
$Na$ as $Na$ is normal over $A$) and $b_0 \in M \setminus Na$. Thus by Lemma
\ref {makkai} there is an automorphism of $M$ fixing $A'$ and taking $a$ to
$b_0$. Thus we can assume that $b_0 \in A$. Iterating this construction we can
assume that $B \subseteq A$.

But now the claim follows from the Open Mapping Theorem. Indeed if $c \in
\cl_{\mathfrak p}(B) \cap M$, then $A \nforks_B c$ and hence $c \nforks_B A$.
Since $\tp(c/A)$ is isolated so is $\tp(c/B)$.
\end {proof}

Now we are ready to prove the main result of this section.

\begin {theorem}
\label {ostable-main}
Let $\mathcal C$ be the class of prime models over Morley sequences in
$\mathfrak p$. Then $\mathcal C$ together with restrictions of $\cl_{\mathfrak
p}$ is a quasiminimal excellent class.
\end {theorem}

\begin {proof}
The only axiom not covered so far is the $\aleph_0$-homogeneity over countable
closed models. Let $H, H' \in \mathcal C$ be prime over Morley sequences $A$
and $A'$ respectively. Let $G \subseteq H, G' \subseteq H'$ be countable closed
or empty. Let $g : G \to G'$ be an isomorphism. If $G$ (and hence $G'$) is
nonempty, we may assume by Proposition \ref {qpg-prime} that $G =
\cl_{\mathfrak p}(A_0) \cap H$, $G' = \cl_{\mathfrak p}(A_0') \cap H'$ where
$A_0 \subseteq A$, $A_0' \subseteq A$ and $g$ maps $A_0$ to $A_0'$. 

Let $f : H \to H'$ be a partial embedding with a finite domain. Assume that
$\dom(f) = \bar a \bar b$ where $\bar a$ is independent over $G$ and $\bar b
\in \cl_{\mathfrak p}(G\bar a)$. Then $\cl_{\mathfrak p}(G\bar a) \cap H$ is
prime and constructible over $G\bar a$. (In case if $G$ is nonempty,
$\cl_{\mathfrak p}(G\bar a) \cap H$ is prime over $A_0\bar a$ by Proposition
\ref {qpg-prime} and $G\bar a$ is normal over $A_0\bar a$.) Since $\tp(\bar
b/G\bar a)$ is isolated, $\cl_{\mathfrak p}(G\bar a) \cap H$ is also
constructible and so prime over $G\bar a\bar b$. Thus the elementary mapping $g
\cup f$ extends to any element of $\cl_{\mathfrak p}(G\bar a\bar b) \cap H$.
\end {proof}

\section {Conclusion}

We have shown that $\beth_{\omega_1}$ is an upper bound for the Hanf number of
the existence of quasiminimal models. We don't have an example to show that it
is sharp, we don't even know if the Hanf number is $\aleph_2$. One way of
proving it would be to show that the existence of a symmetric strongly regular
type is enough to construct arbitrarily large quasiminimal models. We have shown
so using additional assumptions: either stability or presence of definable
Skolem functions. One could try to expand the language with Skolem functions,
but doing so naively will destroy the regular type. Detailed examination of the
proof however reveals that not all Skolem functions are necessary, which gives
this approach some hope.

On the other direction one could ask whether the stability assumptions we use
for our constructions are sharp. For Theorem \ref {stable-main} one can ask
whether there is a simple theory with a strongly regular type and no arbitrarily
large quasiminimal model. A candidate for such a theory could be ACFA where the
type of a transformally transcendental element is strongly regular. By Theorem
\ref {gen-main} there is a quasiminimal model of ACFA of cardinality $\aleph_1$,
but we don't know what happens in other cardinalities.

One can also try to construct a quasiminimal excellent class from a weaker
assumption than $\omega$-stability. A reasonable assumption could be a strongly
regular type in a superstable theory with NOTOP. There is a deep analogy between
NOTOP and excellence and one could try to exploit it. The methods of \cite
{bhhkk} allow one to construct a quasiminimal excellent class if the theory is
superstable and all types are definable over finite subsets. This is however not
very far from an $\omega$-stable theory.  Every such theory which is in addition
small (i.e. has countably many pure types) is $\omega$-stable.

\bibliographystyle {plainnat}
\bibliography {../all}
\end {document}